\documentclass[12pt]{article}
\usepackage{amssymb}
\usepackage{amsmath}
\newtheorem{prop}{Proposition}[section]
\newtheorem{thm}{Theorem}

\newtheorem{remark}[prop]{Remark}

\newtheorem{cor}[prop]{Corollary}
\newtheorem{ex}[prop]{Example}

\newcommand{\Affine}{\mbox{\rm Affine}}
\newcommand{\fix}{\mbox{\rm Fix}}

\def\R{\Bbb R} \def\N{\Bbb N}  \def\Z{\Bbb Z} 
 
  \def\Q{\Bbb Q}  

\begin{document}
\title{Free affine $\Z^p$-actions on Tori.}

\author{Richard Urz\'{u}a-Luz \footnote{ Partially
supported by Fondecyt \# 1100832 and DGIP of the Universidad Cat\'olica del
Norte, Antofagasta, Chile}}

\date{}
\maketitle
\begin{abstract}

We prove that any $\Z^p$-action ${\bf A}$ that acts by automorphisms of  $\Z^q$ with a non-zero fixed-point set induces a unipotent factor of the  $\Z^p$-action ${\bf A}$ which determines whether the action ${\bf A}$ is {\it liberated affine}, i.e. ${\bf A}$ is the linear part of a free affine $\Z^p$-action on the torus $T^q$.
 In general, it is not true that all unipotent $\Z^p$-actions $\bf{U}$ on $\Z^q$ are liberated affine: counter-examples appear for $q\geq 4$.
But if the dimension of the  fixed-point set of ${\bf U}$ regarded as a subspace of $\Z^q$ is less than $q/2$, then ${\bf U}$ is liberated affine.
If $q\leq3$, then a $\Z^p$-action on $\Z^q$ with non-zero fixed-point set is liberated affine.
Finally, for unipotent $\Z^p$-actions on $\Z^3$, we obtain a classification of all those that are the linear part of a minimal free affine $\Z^p$-action on $T^3$.

\end{abstract}
\section{\textbf{Introduction.}}

 Let $R$ be a ring with unity. Denote by Hom$_R(R^{q_1},R^{q_2})$  the $R$-module of all homomorphisms of $R^{q_1}$ to $R^{q_2}$.
It is known that any $\Z^p$-action ${\bf A}$ that acts by automorphisms of $\Z^q$  with non-trivial fixed-point set decomposes $\Z^q$ into a direct product $\Z^q=\Z^{q_1}\times\Z^{q_2}$ such that
 \begin{equation}
 {\bf A}(\ell)(x,y)=({\bf A }_1(\ell)(x),{\bf V}(\ell)(x)+{\bf A}_2(\ell)(y))
 \label{I1}
 \end{equation}
 where ${\bf A}_1$ is an unipotent $\Z^p$-action on $\Z^{q_1}$, ${\bf A}_2$ is a $\Z^p$-action on $\Z^{q_2}$  acting without fixed points on $\Z^{q_2}-\{0\}$, and ${\bf V}$ is a map of $\Z^p$ to  Hom$_\Z(\Z^{q_1},\Z^{q_2})$ that satisfies the following 1-cocycle property.
 \begin{equation}
 {\bf V}(\ell+\ell')={\bf V}(\ell)\circ{\bf A}_1(\ell')+{\bf A}_2(\ell)\circ{\bf V}(\ell')
 \end{equation}
 for all $\ell$ and $\ell'$ in $\Z^p$.
 The $\Z^p$-actions ${\bf A}_1$ and ${\bf A}_2$  induce a $\Z^p$-bimodule structure on the vector space Hom$_\Q(\Q^{q_1},\Q^{q_2})$ given by composition on the right and by composition on the left respectively.
The $\Z$-module Hom$_\Z(\Z^{q_1},\Z^{q_2})$ has a $\Z^p$-bimodule structure as a subbimodule of Hom$_\Q(\Q^{q_1},\Q^{q_2})$.

In Section 3 we prove the following.
\begin{thm}\label{reducida}
Let ${\bf A }$ and ${\bf A}_1$  be the $\Z^{p}$-actions defined above.
Then, the action ${\bf A}$ is liberated affine if, and only if, the unipotent action ${\bf A}_1$ is liberated affine.
\end{thm}
Suppose that the 1-cocycle ${\bf V}$ is trivial in the 1-cohomology of Hom$_\Z(\Z^{q_1},\Z^{q_2})$. Then there is a $W_0\in \mbox{Hom}_\Z(\Z^{q_1},\Z^{q_2})$ such that
\begin{equation}
{\bf V}(\ell)=W_0\circ{\bf A}_1(\ell)-{\bf A}_2\circ(\ell)W_0 \mbox{ for all }\ell\in\Z^p.
\label{I2}
\end{equation}
Then the action ${\bf A}$ given in (\ref{I1}) is conjugated by the automorphism of $\Z^q$, $H(x,y)=(x,-W_0(x)+y)$ to the $\Z^p$-action ${\bf B}$ on $\Z^q=\Z^{q_1}\times\Z^{q_2}$ given by
\begin{equation}
{\bf B}(\ell)(x,y)=({\bf A}_1(\ell)(x),{\bf A}_2(\ell)(y))
\end{equation}
in which case Theorem 1 is immediate. Indeed, since theorem 4 of Hirsch says that all 1-cocycles over the action ${\bf A}_2$ with coefficients in $\R^{q_2}$ are coboundaries, then all affine $\Z^p$-actions $\varphi$ on the torus with linear part ${\bf B}$ are conjugate by a translation to
\begin{equation}
\varphi_0(\ell)(x,y)=({\bf A}_1(\ell)(x)+\alpha(\ell), {\bf A}_2(\ell)(y)).
\end{equation}
However, corollary \ref{bimodulo} in appendix says that (\ref{I2}) has a solution in Hom$_\Q(\Q^{q_1},\Q^{q_2})$ and not  in Hom$_\Z(\Z^{q_1},\Z^{q_2})$, and so the automorphism $H$ does not always induce an automorphism of  $\Z^q$.
Therefore, the previous proof does not consider all possible cases.
In the proof  of Theorem 1 we use corollary 6.1  and the theorem 4 of Hirsch,  specifically that all the 1-cocycles and 2-cocycles over the action ${\bf A}_2$ with coefficients in $\R^{q_2}$ are coboundaries.

In Section 2 we make some definitions and provide some general results about the affine actions on tori, and we show that to study of free affine actions (minimal free affine actions or 1-cohomologically rigid minimal free actions) on the torus, we can assume, without loss of generality,  that our affine actions on the torus are isotopic to their linear part.
This considerably simplifies the proofs of the results and the construction of the examples in Sections 3, 4 and 5. However, for the non-affine case, the isotopy class of the action
plays an important role in the construction of a free $\Z^p$-action  on the torus $T^3$ whose linear part acts without fixed points on $\Z^3-\{0\}$. This example in ~\cite{L1} gives a negative answer to the {\it strong question} posed in ~\cite[Problem 2.8]{S}.

 %In ~\cite{L1} we construct an example of a free $\Z^p$-action $\varphi$  on the torus $T^3$ such that $\fix (\varphi_\ast)= \{0\}$ and that is not isotopic to its linear part. We do not know if there is an example that is isotopic to its linear part. In general, isotopy class will be important in the study of the actions on tori.

In Section 4, Example \ref{ejemplo} shows that, in general, an unipotent $\Z^p$-action ${\bf U}$ on $\Z^q$ is not necessarily liberated affine.
These examples begin to appear for $p\geqslant 3$ and $q \geqslant 4 $.
Proposition 3.3 shows that Example \ref{ejemplo} does not exist for $p = 2$.
 On the other hand,  if $2\dim \fix({\bf U})>q$, then the unipotent $\Z^p$-action ${\bf U}$ on $\Z^q$ is liberated affine.
We do not know whether any unipotent $\Z^2$-action on $\Z^q$ is liberated affine.

In Section 5 we show that, for low dimension, any  $\Z^p$-action on $\Z^q$ with non-trivial fixed-point set is liberated affine.
This result complements ~\cite[Theorem 1]{L1}].
In subsection 5.1 we obtain a classification of all those actions that are the linear part of a minimal free affine $\Z^p$-action on the three dimensional torus.
This result solves Problem 1 proposed in ~\cite{L} for low-dimensional tori.

\section{\textbf{Preliminaries}}

Consider the $q$-dimensional torus $T^q=\R^q/\Z^q$ and let $\varphi$ be a $\Z^{p}$-action on $T^q$.
We denote by $\varphi_\ast$ the $\Z^p$-action induced by $\varphi$ on the first homology group $H_1(T^q,\Z)=\Z^q$.
 The $\Z^{p}$-action $\varphi$ is said to be {\it free} if $\varphi (\ell)$ has no fixed points, for each  $\ell \in \Z^p-\{0\}$, it follows from the Lefchetz fixed-point theorem that $1$ is an eigenvalue of $\varphi_*(\ell)$ for each $\ell\in\Z^p$.
We denote the set of fixed points of $\varphi_*(\ell)$ by $\fix (\varphi_\ast)$ , i.e.
\begin{equation}
\fix ({\varphi_\ast})=\{k\in \Z^{q}|{\bf \varphi_\ast}(\ell )k=k \text{ for all }\ell \in \Z^{p}\}.
\end{equation}
Let $\Affine (T^q)$ denote the group of affine transformations of the torus $T^q$.
By an affine  $\Z^{p}$-action on $T^q$
we mean a homomorphism $\varphi$ of $\Z^p$ into $\Affine (T^q)$.
In~\cite{S}, dos Santos proved that every affine $\Z^{p}$-action $\varphi$ on the torus $T^{q}$, $q\in\N$ satisfying $\fix (\varphi_*)=\{0\}$ has a finite orbit.
In particular, $\varphi$ is not free.
Hence all free affine $\Z^p$-actions $\varphi$ on the torus $T^q$ induce a $\Z^{p}$-action  $\varphi_\ast$ on $H_{1}(T^{q},\Z)$ such that $\fix (\varphi_\ast)\neq \{0\}$.
 A $\Z^p$-action ${\bf A}$ that acts by automorphisms of $\Z^q$ is said to be {\it liberated affine}
if there is a free affine $\Z^{p}$-action $\varphi$ on the torus $T^{q}$ such that $\varphi_\ast={\bf A}$.
\begin{remark}
In this paper all the $\Z^p$-actions on $\Z^q$ act by automorphisms of $\Z^q$.
\end{remark}

Let ${\bf A}$ be a $\Z^p$-action on $\Z^q$.
A map $\alpha:\Z^p\rightarrow T^q$ is a 1-cocycle over the action ${\bf A}$ if for each $\ell,\ell'\in\Z^p$
\begin{equation}
\alpha(\ell+\ell')={\bf \alpha}(\ell){\bf +}{\bf A}(\ell)({\bf \alpha}(\ell'))
\end{equation}
A 1-cocycle is trivial if there exists an ${\bf x}_0\in T^q$ such that
\begin{equation}
\alpha(\ell)={\bf x}_0{\bf -}{\bf A}(\ell)({\bf x}_0)
\end{equation}
where the automorphisms ${\bf A}(\ell)$ of $\Z^q$ can be interpreted as automorphisms of the torus $T^q$ for all $\ell\in\Z^p$.

The 1-cohomology group $H^1(\Z^p,T^q)$ is the quotient space of the 1-cocycles modulo the trivial 1-cocycles.
It is well known that the basic structure of affine actions can be expressed in terms of the 1-cohomology.
Given an affine $\Z^p$-action $\varphi$ on the torus with linear part ${\bf A}$, define the 1-cocycle by its action on the identity element, $\alpha(\ell)=\varphi(\ell)(0)$, for $\ell\in\Z^p$.
Conversely.
given a 1-cocycle $\alpha:\Z^p\rightarrow T^q$ over the action {\bf A}, define an affine  $\Z^p$-action by the rule $\varphi(\ell)({\bf x})=\alpha(\ell)+{\bf A}(\ell)({\bf x})$ for ${\bf x}\in T^q$ and $\ell\in\Z^p$.
Let $\alpha$ be the 1-cocycle over the action ${\bf A}$ associated to an action $\varphi$.
The 1-cocycle class associated to $\varphi$ is the cohomology class $[\alpha]\in H^1(\Z^p,T^q)$.
The existence of a fixed point for an affine action on a torus is related to a group extension problem and the vanishing of the corresponding cohomology groups: for example, the class $[\alpha]$ is trivial if, and only if, the action $\varphi$ has a fixed point; we obtain a dense set of periodic points if, and only if, $[\alpha]$ is torsion. For more information, see ~\cite{Hu}.

In the study of the free affine $\Z^p$-actions on a torus $T^q$ we consider the long exact sequence of cohomology groups associated with the exact sequence of $\Z^p$-modules
\begin{equation}
0\rightarrow \Z^q\rightarrow\R^q\rightarrow T^q\rightarrow 0,
\end{equation}

\begin{equation}
\cdot\cdot\cdot\rightarrow H^1(\Z^p,\R_{{\bf A}}^q)\overset{\pi}{\rightarrow} H^1(\Z^p,T_{{\bf A}}^q)\overset{\delta}{\rightarrow} H^2(\Z^p,\Z_{{\bf A}}^q)\overset{i}{\rightarrow} H^2(\Z^p,T_{{\bf A}}^q)\rightarrow\cdot\cdot\cdot
\end{equation}

By ~\cite{Hu}, we can identify the kernel of  $H^2(\Z^p,\Z_{{\bf A}}^q)\overset{i}{\rightarrow} H^2(\Z^p,\R_{{\bf A}}^q)$ with the torsion subgroup of  $H^2(\Z^p,\Z_{{\bf A}}^q)$.
Suppose given an affine $\Z^p$-action $\varphi$ on $T^q$ with linear part ${\bf A}$ and consider, for each $\ell\in\Z^p$, a lifting $\widetilde{\varphi}_1(\ell):\R^q\rightarrow\R^q$ of $\varphi(\ell):T^q\rightarrow T^q$ to the covering $\R^q$.
Then $\widetilde{\varphi}_1(\ell)$ can be written as
\begin{equation}
\widetilde{\varphi}_1(\ell)(x)={\bf A}(\ell)(x)+\widetilde{\alpha}(\ell)
\end{equation}
where $x\in\R^q$.
The mapping $\widetilde{\alpha}:\Z^p\rightarrow \R^q$ is not in general a cocycle over the action ${\bf A}$, but as it comes from the 1-cocycle $\varphi(\ell)(0), \ell\in\Z^p$, we must have
\begin{equation}
\widetilde{\alpha}(\ell+\ell')=\widetilde{\alpha}(\ell)+{\bf A}(\ell)(\widetilde{\alpha}(\ell'))-k(\ell,\ell')
\end{equation}
where $k:\Z^q\times\Z^q\rightarrow\Z^q$.
Since $k(\ell,\ell')=\widetilde{\alpha}(\ell)-\widetilde{\alpha}(\ell+\ell')+$ ${\bf A}(\ell)\widetilde{\alpha}(\ell')$ for all $\ell,\ell'\in\Z^p$, it follows that $k$ is a trivial 2-cocycle in $H^2(\Z^p,\R_{{\bf A}}^q)$ and so $k$ is torsion in $H^2(\Z^p,\Z_{{\bf A}}^q)$, i.e. there are $s\in\N$ and $\widetilde{\beta}:\Z^p\rightarrow \Z^q$ such that
\begin{equation}
sk(\ell,\ell')=\widetilde{\beta}(\ell)-\widetilde{\beta}(\ell+\ell')+{\bf A}(\ell)(\widetilde{\beta}(\ell'))
\end{equation}
for all $\ell,\ell'\in\Z^p$.
Now consider the affine $\Z^p$-action $\psi$ on torus $T^q$ given on the covering $\R^q$ by
\begin{equation}
\widetilde{\psi}(\ell)(x)={\bf A}(\ell)(x)+\widetilde{\alpha}_1(\ell)
\label{fix}
\end{equation}
where
\begin{equation}
\widetilde{\alpha}_1(\ell)=\widetilde{\alpha}(\ell)-\frac{1}{s}\widetilde{\beta}(\ell),\mbox{  }\ell\in\Z^p.
\label{fix1}
 \end{equation}
 By construction, $\widetilde{\alpha_1}:\Z^p\rightarrow \R^q$ is a 1-cocycle over the action ${\bf A}$, therefore $\widetilde{\psi}$ is an affine $\Z^p$-action on $\R^q$.
Therefore, the affine $\Z^p$-action $\psi$ is isotopic to its linear part ~\cite{S}.
We call $\psi$ the {\it principal affine $\Z^p$-action} of $\varphi$.
 \begin{remark}
 Let $\Gamma$ be a finitely generated group. The principal affine action can be defined in the same way for all affine $\Gamma$-actions on tori.
 \end{remark}

 The following proposition shows that if $\varphi$ is a free affine $\Z^p$-action or minimal $\Z^p$-action or 1-cohomologically rigid $\Z^p$-action, then its principal affine $\Z^p$-action $\psi$ also has these properties.
Of course the converse is also true.
\begin{prop}\label{principal}
Let $\varphi$ be an affine $\Z^p$-action on $T^q$ and $\psi$ the  principal affine $\Z^p$-action of $\varphi$. Then:
\begin{enumerate}
\item [i.] If $\varphi$ is a free $\Z^p$-action, then $\psi$ is a free $\Z^p$-action.
\item [ii.] If $\varphi$ is a minimal $\Z^p$-action, then $\psi$ is a minimal $\Z^p$-action.
\item[iii.] If $\varphi$ is a 1-cohomologically rigid $\Z^p$-action, then $\psi$ is too.
\end{enumerate}
\end{prop}
\noindent{\bf Proof of (i).} Suppose that $\psi$ is a $\Z^p$-action that is not free. Then there are ${\bf x}_0\in T^q$ and $\ell_0\in\Z^p-\{0\}$ such that
\begin{equation}
\psi(\ell_0)({\bf x}_0)={\bf x}_0.
\label{fix3}
\end{equation}

On the covering $\R^q$ $\psi(\ell_0)$,  we can write $\widetilde{\psi}(\ell_0)={\bf A}(\ell_0)+\widetilde{\alpha}_1(\ell)$. Then from Eq.~(\ref{fix}) we have that
\begin{equation}
\widetilde{\psi}(\ell_0)(\widetilde{x}_0)={\bf A}(\ell_0)\widetilde{x}_0+\widetilde{\alpha}_1(\ell_0)=\widetilde{x}_0+N
\label{fix2}
\end{equation}
where $\widetilde{x}_0$ is a lifting of ${\bf x}_0$ to $\R^q$ y $N\in\Z^q$.
By (\ref{fix3}) and (\ref{fix2}), we get
\begin{equation}
{\bf A}(\ell_0)(\widetilde{ x}_0)+\widetilde{\alpha}(\ell_0)= \widetilde{ x}_0+\frac{1}{s}\widetilde{\beta}(\ell_0)+N.
\end{equation}
I.e. $\widetilde{\varphi}(\ell_0)(\widetilde{x}_0)=\widetilde{x}_0+\frac{1}{s}\widetilde{\beta}(\ell_0)+N$.
Hence $\widetilde{\varphi}^n(\ell_0)(\widetilde{x}_0)\in \widetilde{x}_0+\frac{1}{s}\Z^q$ for all $n\in\Z$, and so the orbit of ${\bf x}_0$ under $\varphi(\ell_0)$ is contained in a finite subset of $T^q$, meaning $\varphi$  is not free.
\newline
\noindent{\bf Proof of (ii).} The proof  follows from ~\cite[Theorem 1]{L} that  $\varphi$ satisfies the irrationality condition on $\Gamma$ if, and only if, $\psi$ does.\newline
\noindent{\bf Proof of (iii).}  It is not difficult to see that $\varphi$ satisfies the hypothesis of \cite[Theorem 2]{L} if, and only if, $\psi$ does.
$\hfill \square$
\begin{remark}
For convenience, in this paper, we assume that any affine action is isotopic to its linear part.
\end{remark}

\section{\textbf{Free affine $\Z^p$ actions on tori and unipotent $\Z^p$-actions}}

The main goal of this section is to show that the study of the free affine $\Z^p$-actions  on a torus be reduced to studying the free affine $\Z^p$-actions on that torus whose linear part is a unipotent action.
 Let ${\bf A}$ be  a $\Z^p$-action on $\Z^q$ that acts by automorphisms of $\Z^q$ and whose fixed-point set $\fix({\bf A})\neq\{0\}$.
It is known that we can decompose $\Z^q$ into a direct product $\Z^q=\Z^{q_1}\times\Z^{q_2}$ such that
 \begin{equation}
 {\bf A}(\ell)(x,y)=({\bf A }_1(\ell)(x),{\bf V}(\ell)(x)+{\bf A}_2(\ell)(y))
 \label{aff1}
 \end{equation}
 where ${\bf A}_1$ is an unipotent $\Z^p$-action on $\Z^{q_1}$, ${\bf A}_2$ is a $\Z^p$-action on $\Z^{q_2}$  acting without fixed points on $\Z^{q_2}-\{0\}$, and ${\bf V}$ is a map of $\Z^p$ to Hom$_\Z(\Z^{q_1},\Z^{q_2})$ that satisfies the following 1-cocycle property.

 \begin{equation}
 {\bf V}(\ell+\ell')={\bf V}(\ell)\circ{\bf A}_1(\ell')+{\bf A}_2(\ell)\circ{\bf V}(\ell')
 \label{aff2}
 \end{equation}
 for all $\ell$ and $\ell'$ in $\Z^p$.

Corollary \ref{bimodulo}  of the Appendix shows that the 1-cocycle ${\bf V}$  defined above is a trivial 1-cocycle in the $\Z^p$-bimodule Hom$_\Q(\Q^{q_1},\Q^{q_2})$, that is, there is a $W_0\in\mbox{Hom}_\Q(\Q^{q_1}.\Q^{q_2})$ such that ${\bf V}(\ell)=W_0\circ{\bf A}_1(\ell)-{\bf A}_2(\ell)\circ W_0$

\noindent{\bf Proof of Theorem \ref{reducida}.} Let us suppose that $\varphi$ is an affine $\Z^p$-action on $T^q$ whose induced action on $H_1(T^q, \Z)$ is ${\bf A}$ and has a lifting to an affine $\Z^p$-action $\widetilde{\varphi}$ on the covering $\R^q$ given by
\begin{eqnarray}
\widetilde{\varphi}(\ell)(x,y)&=&(\widetilde{\varphi}_1(\ell)(x), {\bf V}(\ell)(x)+{\bf A}_2(\ell)(y)+\beta(\ell))
\\
&=&({\bf A}_1(\ell)(x)+\alpha(\ell), {\bf V}(\ell)x+{\bf A}_2(\ell)(y)+\beta(\ell))
\label{33c}
\end{eqnarray}
where $x\in \R^{q_1}$, $y\in \R^{q_2}$ and $\ell\in\Z^p$.
Since $\widetilde{\varphi}$ is isotopic to its linear part, we can assume that the map $\alpha$ of $\Z^p$ to  $\R^{q_1}$ is a 1-cocycle  over the action ${\bf A}_1$  and $\beta$ is a map of $\Z^p$ to $\R^{q_2}$ satisfying
\begin{equation}
\beta(\ell+\ell')={\bf A}_2(\ell)(\beta(\ell'))+\beta(\ell)+{\bf V}(\ell)(\alpha(\ell')) \mbox{ for all } \ell, \ell'\in\Z^p.
\label{4c}
\end{equation}
By Corollary \ref{bimodulo},  we show that $\beta(\ell)-W_0(\alpha(\ell)),\ell\in\Z^p$ is a 1-cocycle over the action ${\bf A}_2$ and so, by  Theorem \ref{hirsch} in the Appendix, is a trivial 1-cocycle.
Thus there is a $z_0\in \R^{q_2}$ such that  $\beta(\ell)-W_0(\alpha(\ell))=(I-{\bf A}_2(\ell))(z_0)$ for all $\ell\in\Z^p$.
Therefore, we can write (\ref{33c})  as
$$
\widetilde{\varphi}(\ell)(x, y)=({\bf A}_1(\ell)(x)+\alpha(\ell),
$$
\begin{equation}
{\bf V}(\ell)(x)+{\bf A}_2(\ell)(y)+W_0(\alpha(\ell))+(I-{\bf A}_2(\ell))(z_0))
\label{34c}
\end{equation}

Now, let $\varphi_1$ be the affine $\Z^p$-action on $T^{q_1}$ given on the covering $\R^{q_1}$ by $\widetilde{\varphi_1}(\ell)={\bf A}_1(\ell)+\alpha(\ell)$.
Suppose  that  $\varphi_1$ is an action that is not free. Then on the covering $\R^{q_1}$ there are $x_0\in \R^{q_1}$, $\ell_0\in\Z^p-\{0\}$ and $n_0\in\Z^{q_1}$ such that
\begin{equation}
\widetilde{\varphi_1}(\ell_0)(x_0)={\bf A}_1(\ell_0)(x_0)+\alpha(\ell_0)=x_0+n_0
\label{5c}
\end{equation}

Putting $\ell=\ell_0$ and $x=x_0$ in (\ref{34c}) we have that
\begin{equation}
\widetilde{\varphi}(\ell_0)(x_0, y)=(x_0+n_0, {\bf V}(\ell_0)(x_0)+{\bf A}_2(\ell_0)(y)+W_0(\alpha(\ell_0))+(I-{\bf A}_2(\ell_0))(z_0))
\end{equation}
and again by Corollary \ref{bimodulo} and Eq.~(\ref{5c}) we obtain
\begin{equation}
\widetilde{\varphi}(\ell)(x_0,y)=(x_0+n_0, A_2(\ell_0)(y)+ (I-{\bf A}_2(\ell_0))(z_0+W_0x_0)+W_0(n_0))
\end{equation}
Now writing $y_0$ for the point $z_0+W_0(x_0)$, we have that $\widetilde{\varphi}(\ell_0)(x_0,y_0)=(x_0+n_0, y_0+W_0(n_0))$, hence the orbit $\{\varphi^m(\ell_0)(x_0,y_0)/m\in\Z\}$ is finite, and so $\varphi$ is not free.

 Conversely, let us suppose that $\varphi_1$ is a free affine $\Z^{p}$-action on $T^{q_1}$ whose induced action on $H_1(T^{q_1}, \Z)$ is ${\bf A}_1$.

Let $\varphi$ be the affine $\Z^p$-action on $T^{q}$ given on the covering $\R^{q}$ by (\ref{33c}) with $\widetilde{\varphi}_1(\ell)={\bf A}_1(\ell)+\alpha(\ell)$.
Since $\varphi_1$ is isotopic to its linear part, we can assume that $\alpha$ is a 1-cocycle over the action ${\bf A}_1$.

To show that  $\varphi$ is a free affine $\Z^p$-action, it is enough to prove that Eq.~(\ref{4c})
admits a solution $\beta:\Z^p\rightarrow \R^q$, and since  $\varphi_1$ is a free $\Z^p$-action, then $\varphi$ is a free $\Z^p$-action.

Indeed, the map $\rho:\Z^p\times\Z^p\rightarrow\R^{q_2}$ defined by $\rho(\ell,\ell')=-{\bf V}(\ell)(\alpha(\ell'))$, for all $\ell,\ell'\in\Z^p$, is a
 2-cocycle over the action ${\bf A}_2$ since
$$
 {\bf A}_2(\ell_1)(\rho(\ell_2,\ell_3))-\rho(\ell_1+\ell_2,\ell_3)+\rho(\ell_1,\ell_2+\ell_3)-\rho(\ell_1,\ell_2)=0.
 $$
The equation above follows directly  from Eqs~(\ref{aff2}) and the fact that $\alpha$ is a 1-cocycle over the action ${\bf A}_1$.
Since $ \fix ({\bf A}_2)=\{0\}$ and by  Theorem \ref{hirsch} of the Appendix,  we have that $H^2(\Z^p;\R^{q_2}_{{\bf A}_2})$ is trivial, thus $\rho$ is a coboundary, i.e. there is a $\beta:\Z^p\rightarrow \R^{q_2}$  satisfying  (\ref{4c}).

$\hfill \square$

\begin{cor}\label{subreducida}
Let $\mathbf{A}$ be a $\Z^p$-action on $\Z^{q}$ with $\fix \left( \mathbf{A}\right) \neq \{0\}.$ Assume that ${\bf A}_1$  is the identity $I$ in Theorem \ref{reducida}.
Then the action $\mathbf{A}$ is liberated affine.
\end{cor}
\noindent{\bf Proof.}
This corollary is immediate because the action ${\bf A}_1(\ell)=I$ for all $\ell\in\Z^p$ is liberated by translations, i.e. defined by the $\Z^p$-action $\varphi$ on $T^q$ given on the covering $\R^{q_1}$ by  $\varphi_1(\ell)=I+\alpha(\ell)$ where $\alpha$ is a linear homomorphism  of $\Z^p$ in $\R^{q_1}$ such that for all $\ell\in\Z^p-\{0\}$, $<\alpha(\ell), k>\notin \Z$ for some $k\in\Z^{q_1}-\{0\}$.
$\hfill \square$
\section{\textbf{ Unipotent Actions}}

The example below shows that in general,  not every unipotent $\Z^p$-action ${\bf U}$ on $\Z^q$ is liberated affine.
\begin{ex}\label{ejemplo}
Let ${\bf U}$ be the unipotent $\Z^p$-action on  $\Z^{2n}=\Z^n\times\Z^n$  defined by
\begin{equation}
{\bf U}(\ell)(x,y)=(x,y+{\bf V}(\ell)(x))
\end{equation}
where ${\bf V}$ is a map of $\Z^p$ to $\mbox{Hom}_\Z(\Z^n,\Z^n)$ that satisfies the following 1-cocycle property.
 \begin{equation}
 {\bf V}(\ell+\ell')={\bf V}(\ell)+{\bf V}(\ell')
 \label{aff3}
 \end{equation}
for all $\ell,\ell'\in\Z^p$.
Consider the affine $\Z^p$-action on $T^{2n}$ given on the covering $\R^{2n}=\R^n\times\R^n$ by
$$
\widetilde{\varphi}(\ell)(x,y)={\bf U}(\ell)(x,y)+\mathbf{\delta}(\ell)
=(x+\alpha(\ell),{\bf V}(\ell)(x)+y+\beta(\ell))\;\; \ell\in\Z^p.
$$
where $\mathbf{\delta}(\ell)=(\alpha(\ell),\beta(\ell))$.
Since $\widetilde{\varphi}$ is an affine $\Z^p$-action on $\R^{2n}$,
\begin{equation}
\alpha(\ell+\ell')=\alpha(\ell)+\alpha(\ell')
\end{equation}
and
\begin{equation}
{\bf V}(\ell)(\alpha(\ell'))={\bf V}(\ell')(\alpha(\ell))
\label{6c}
\end{equation}
for all $\ell, \ell'\in\Z^p$.
Suppose there is an $\ell_0\in\Z$ such that ${\bf V}(\ell_0)$ is invertible.
Then,  (\ref{6c}) implies that
\begin{equation}
\alpha(\ell)={\bf V}(\ell_0)^{-1}\circ{\bf V}(\ell)(\alpha(\ell_0))
\label{7c}
\end{equation}
for all $\ell\in\Z^p$.
By \ref{6c} and \ref{7c} we have
\begin{equation}
{\bf V}(\ell_1)\circ{\bf V}(\ell_0)^{-1}\circ{\bf V}(\ell_2)(\alpha(\ell_0))={\bf V}(\ell_2)\circ{\bf V}(\ell_0)^{-1}\circ{\bf V}(\ell_1)(\alpha(\ell_0))
\end{equation}
Defining the operation ${\bf V}(\ell)\star{\bf V}(\ell')={\bf V}(\ell)\circ{\bf V}(\ell_0)^{-1}\circ{\bf V}(\ell')$ and choosing the commutator $[{\bf V}(\ell_1),{\bf V}(\ell_2)]$ to be invertible, we obtain that $\alpha(\ell_0)=0$, and a direct calculation shows that the points
$(-{\bf V}(\ell_0)^{-1}(\beta(\ell_0)),y)$ are fixed points of $\varphi(\ell_0)$.
Note that it is necessary that $p\geq3$ and $n\geq2$.
Proposition \ref{rango2} shows that for $p = 2$ there are no such examples.
\end{ex}

Note that all unipotent $\Z^p$-actions $\mathbf{U}$ on $\Z^q$ decompose $\Z^q$ into a direct product $\Z^q=\Z^{q-k}\times\Z^k$ such that
 \begin{equation}
 {\bf U}(\ell)(x,y)=({\bf U }_1(\ell)(x),y+{\bf V}(\ell)(x))
 \label{uni4,5}
 \end{equation}

\noindent where $k=\dim \fix({\bf U})$,  ${\bf U}_1$ is an unipotent $\Z^p$-action on $\Z^{q-k}$, and ${\bf V}$ is a 1-cocycle over the action ${\bf U}_1$.
\begin{thm} \label{rango}
Let ${\bf U}$ be an unipotent $\Z^p$-action on $\Z^q$.
Suppose that the 1-cocycle ${\bf V}$ over the action ${\bf U}_1$ has rank less than $k$, i.e. the linear map
${\bf V}(\ell):\R^{q-k}\rightarrow \R^k$ has rank less than $k$ for all $\ell\in\Z^p$.  Then ${\bf U}$ is liberated affine.
\end{thm}
\noindent{\bf Proof.}

Consider the affine  $\Z^p$-action $\varphi$ on $T^q$ on the covering $\R^q$ by

\begin{equation}
\tilde{\varphi}\left( \ell \right) ={\bf U}\left( \ell \right) +\gamma
\left( \ell \right)
\end{equation}
Take $\gamma \left( \ell \right) =\left(0, \alpha(\ell)\right)$ be where $\alpha :\Z^{p}\rightarrow \R^{k}$ is
an homomorphism defined by the product of matrices
$\left( \alpha _{ij}\right) _{k\times p}$ such that
\begin{equation}
\alpha \left( \ell \right) =\left( \alpha _{ij}\right) _{k\times p}\cdot \ell_{p\times 1}
\end{equation}
for all $\ell^t_{p\times 1} =\ell\in \Z^{p}.$ We choose a matrix $\left( \alpha _{ij}\right) _{k\times p}$ such that the family of real numbers
 $\{\alpha _{ij}; \ 1\leq i \leq k\mbox{ and } 1\leq j \leq p  \}$ are linearly independent over the rationals, in particular, the inner product
 \begin{equation}
 \langle \mathbf{m},\alpha(\ell)\rangle \not \in \Z \mbox{ for all } \mathbf{m}\in \Z^k-\{0\} \mbox{ and for all } \ell\in\Z^p-\{0\}.
 \label{uni6}
 \end{equation}
Now, suppose there is an $\ell\in\Z^p$ and ${\bf x}\in T^q$ such that $\varphi(\ell)({\bf x})={\bf x}$. Then, on the covering $\R^q$, there is an $N\in\Z^q$ such that
\begin{equation}
\tilde{\varphi}\left( \ell \right) \left( \tilde{{\bf x}}\right) ={\bf U}\left( \ell
\right)(\tilde{{\bf x}})+\gamma \left( \ell \right) =x+N
\label{uni5}
\end{equation}
where $\tilde{{\bf x}}=\left( x,y\right)$ is a lifting of ${\bf x}$ to $\R^{q-k}\times \R^{k}$,  $N=\left(n_{1},n_{2}\right) \in \Z^{q-k}\times\Z^k,$ and so Eq.~(\ref{uni5}) can be written as follows.
$$
\begin{array}{c}
{\bf U}_1\left( \ell \right)(x)=x+n_{1} \\
 y+{\bf V}\left( \ell \right)(x)+\alpha \left( \ell \right)=y+n_{2}
\end{array}
$$
Since the rank of ${\bf V}(\ell)$ is less than $k$, there exists an $\mathbf{m_\ell}\in\Z^k-\{0\}$ such that $\langle\mathbf{m_\ell},{\bf V}(\ell)z\rangle=0$,
for all $z\in \R^{q-k}$.
From the above equation we have
$$
\langle\mathbf{m_\ell},\alpha(\ell)\rangle=\langle\mathbf{m_\ell},n_2\rangle \in \Z.
$$
and so (\ref{uni6}) implies $\ell=0$.
$\hfill \square$
\begin{cor}
Let $\mathbf{U}$ be an unipotent $\Z^p$-action on $\Z^q$. If  $2\dim \fix(\mathbf{U})>q$, then $\mathbf{U}$ is liberated affine.
\end{cor}\label{corrank}
\noindent{\bf Proof.}
Since the linear map ${\bf V}(\ell):\R^{q-k}\rightarrow \R^k$ has rank at most $q-k$ and $q-k<k$, the corollary follows from theorem
 \ref{rango}.

$\hfill \square$

We do not know if all unipotent $\Z^2$-actions on $\Z^q$ are liberated affine.
However, the following
proposition gives a positive response when ${\bf U}_1=I$ for all $\ell$ in Eq.~(\ref{uni4,5}).
\begin{prop}\label{rango2}
Let $\mathbf{U}$ be an an  $\Z^2$-action on $\Z^q$.
Suppose that ${\bf U}_1(\ell)=I$ for all $\ell\in\Z^2$ in Eq.~(\ref{uni4,5}).
Then ${\bf U}$ is liberated affine.
\end{prop}
\noindent{\bf Proof.}

By Theorem \ref{rango} we may assume that there exists  an $e_1\in\Z^2$ such that the linear map ${\bf V}(e_1):\R^{q-k}\rightarrow
\R^k$ has rank $k$.
Let $\alpha(e_1)$ be a vector in $\R^{q-k}$ such that
\begin{equation}
\langle\alpha(e_1),\mathbf{n}\rangle\not \in \Z
\label{7.1}
\end{equation}
for all $\mathbf{n}\in\Z^{q-k}-\{0\}$.
Since for all $\ell\in\Z^2$ there is an $\alpha(\ell)\in\R^{q-k}$ such that
\begin{equation}
{\bf V}(e_1)(\alpha(\ell))={\bf V}(\ell)(\alpha(e_1))
\label{8c}
\end{equation}
we can consider $\alpha(e_2)$ satisfying (\ref{8c}) such that $\{e_1,e_2\}$ is a basis of $\Z^2$.
Define $\alpha(\ell)=l_1\alpha(e_1)+l_2\alpha(e_2)$
where $\ell=l_1e_1+l_2e_2$.
A simple calculation shows that
\begin{equation}
{\bf V}(\ell)(\alpha(\ell'))={\bf V}(\ell')(\alpha(\ell))
\label{9c}
\end{equation}
for all $\ell,\ell'\in\Z^2$.
Now let $\widehat{\beta}:\Z^2\rightarrow\R^k$ be a  homomorphism  such that
\begin{equation}
\widehat{\beta}(\ell)\not\in\Z^k
\label{9.1}
\end{equation}
for all $\ell\in\Z^2-\{0\}$ and define
\begin{equation}
\beta(\ell)={\bf V}(\ell)(\alpha(\ell))+\widehat{\beta}(\ell)
\label{10c}.
\end{equation}

Hence, by construction, putting
\begin{equation}
\widetilde{\varphi}(\ell)(x,y)={\bf U}(\ell)(x,y)+(\alpha(\ell),\beta(\ell))
=(x+\alpha(\ell),y+{\bf V}(\ell)(x)+\beta(\ell))
\label{11c}
\end{equation}
for each  $\ell\in\Z^2$ determines a free affine $\Z^2$-action $\varphi$ on $T^q$.
In fact, suppose we have $\varphi(\ell){\bf x}={\bf x}$ for
some $\ell\in\Z^2$ and some $\mathbf{x}\in T^q$. Then, on the covering $\R^q$, we obtain by (\ref{11c})
\begin{equation}
\widetilde{\varphi }\left( \ell \right) (x,y)=\left( x+\alpha \left( \ell
\right) ,y+{\bf V}\left( \ell \right)(x)+\beta \left( \ell \right) \right) =\left(
x,y\right) +\left( N_{1},N_{2}\right)
\label{12c}
\end{equation}
where $\left( x,y\right) $ is a lifting of ${\bf x}$ to $\R^{q-k}\times \R^{k}$ and $\left( N_{1},N_{2}\right) \in \Z^{q-k}\times \Z^{k}$. Now, from the first coordinate of Eq.~(\ref{12c}), we obtain  $\alpha(\ell)=N_1$ and by (\ref{8c}), $V(\ell)\alpha(e_1)\in\Z^k$. Hence by (\ref{7.1}), necessarily $V(\ell)=0$.
From the second coordinate of Eq.~(\ref{12c}) we have that
$\widehat{\beta}(\ell)=N_2$ and by $(\ref{9.1})$ $\ell=0$.
$\hfill \square$

\section{\textbf{Free affine actions on low-dimensional tori.}}

Example \ref{ejemplo} shows that when $q\geq4$, not every unipotent $\Z^p$-action ${\bf U}$ on $\Z^q$ is the linear part of a
free affine $\Z^p$-action on $T^q$.
In this section we will show any $\Z^p$-action on $\Z^q$, with $1\leq q\leq3$, such that
  $\fix \left( {\bf A}\right) \neq \{0\}$ admits a free affine $\Z^p$-action $\varphi$ on $T^q$ such that $\varphi_\ast={\bf A}$.
\begin{thm}\label{dimensionbaja}
Let $p\geq2$  and $q\leq 3$. Then a $\Z^p$-action ${\bf A}$ on $\Z^q$ such that $\fix({\bf A})\neq\{0\}$ is liberated affine.
\end{thm}
\noindent{\bf Proof.}
For $q=1$,  $\fix \left( \mathbf{A}\right) \neq \{0\}$ implies that ${\bf A}(\ell)=1$ for all $\ell\in\Z^p$. Then it is clear that $\mathbf{A}$ is liberated by translation.

For $q=2$, since $\fix \left( \mathbf{A}\right) \neq \{0\}$, then on $\Z^2$

\begin{equation}
{\bf A}(\ell)(x,y)=(x,\nu(\ell)y+\mu(\ell)x)
\end{equation}
for all $\ell \in \Z^p$,  where $\mu(\ell)\in\Z$ and $\mu (\ell)=1$ or $-1$.
Consider the case $\nu(\ell)=-1$ for some $\ell\in\Z^p$. Then the theorem follows from Corollary \ref{subreducida}.

If $\nu(\ell)=1$ and $\mu(\ell)=0$ for all $\ell\in\Z^p$, then $\mathbf{A}(\ell)=I$ for all  $\ell\in\Z^p$ and $\mathbf{A}$ liberated by translation on $T^2$.

If $\nu(\ell)=1$ for all $\ell\in\Z^p$  and $\mu$  is a non-zero homomorphism,
consider the affine $\Z^p$-action  $\psi(\ell)$ on $T^2$  given on the covering $\R^2$ by
\begin{equation}
\widetilde{\psi}(\ell)(x,y)=(x+\alpha(\ell),\mu(\ell)x+y+\beta(\ell)).
\label{15c}
\end{equation}
Since we can assume that $\widetilde{\psi}$ is a $\Z^p$-action on $\R^2$, then

\begin{enumerate}
\item [i.] $\alpha(\ell+\ell')=\alpha(\ell)+\alpha(\ell')$, for all $\ell,\ell'\in\Z^p$
\item [ii.]$\beta(\ell+\ell')=\mu(\ell)\alpha(\ell')+\beta(\ell)+\beta(\ell')$, for all $\ell,\ell'\in\Z^p$
\end{enumerate}
Now, exchanging $\ell$ and $\ell'$ in item ii above, we have that
 \begin{equation}
 \mu(\ell)\alpha(\ell')=\mu(\ell')\alpha(\ell), \mbox{ for all}  \ell,\ell'\in\Z^p.
 \label{14c}
 \end{equation}
 Note also that items i and  ii  imply that $\widetilde{\psi}$ in $(\ref{15c})$ is an affine $\Z^p$-action on $\R^2$.
 Let $e_1\in\Z^p$ be such that $\Z^p=\Z e_1\oplus \ker\mu$.
Choose $\alpha$ such that $\alpha(e_1)\notin \Q $. Then by (\ref{14c}), $\ker\alpha=\ker\mu$.
On the other hand, now choose $\beta(\ell)=\mu(\ell)\alpha(\ell)+\widehat{\beta}(\ell)$ for all $\ell\in\Z^p$, where $\widehat{\beta}$ is an homomorphism of $\Z^p$ to $\R$ such that $\widehat{\beta}(\ell)\notin \Z$ for all $\ell\in\ker\mu-\{0\}$.
With these choices, the affine  $\Z^p$-action $\psi$ on the torus is a free action. In fact, if $\mathbf{x}\in T^2$ is a fixed point for $\psi(\ell)$, then on the covering $\R^q$ there is an $N=(n_1,n_2)\in\Z^2$ such that
\begin{equation}
\widetilde{\psi}(\ell)(x,y)=(x+\alpha(\ell),\mu(\ell)x+y+\beta(\ell))=(x,y)+(n_1,n_2)
\end{equation}
where $(x,y)$ is a lifting of ${\bf x}$ to $\R^2$.
From the first coordinate, we get $\alpha(\ell)=n_1$, since $\alpha(e_1)\notin \Q $. Now by (\ref{14c}) we have  $\ell\in\ker\mu$. From the second coordinate we get $\beta(\ell)=\widehat{\beta}(\ell)=n_2$, and this is possible only if $\ell = 0$.

For $q=3$, since $\fix \left( \mathbf{A}\right) \neq \{0\}$, it follows that
\begin{equation}
{\bf A}(\ell)(x,Y)=(x,{\bf V}(\ell)(x)+{\bf B}(\ell)(Y))
\end{equation}
where  $\mathbf{B}$ is a $\Z^p$-action on $\Z^2$.
If $\fix\left( \mathbf{B}\right)=\{0\}$, the theorem follows from Corollary \ref{subreducida}.
If $\fix \left( \mathbf{B}\right) \neq \{0\}$, then
\begin{equation}
{\bf A}(\ell)(x,y,z)=(x,\mu(\ell)x+y,\omega(\ell)x+\nu(\ell)y+\tau(\ell)z)
\end{equation}

If now $\tau(\ell)=-1$ for some $\ell\in\Z^p$, the theorem follows from Theorem \ref{reducida} and the case $q=2$.

Finally, the unipotent case,
\begin{equation}
{\bf A}(\ell)(x,y,z)=(x,\mu(\ell)x+y,\omega(\ell)x+\nu(\ell)y+z)
\end{equation}
for all $\ell\in\Z^p$.
Since ${\bf A}$ is a $\Z^p$-action on $\Z^3$, $\mu$ and $\nu$ are homomorphisms of $\Z^p$ to $\Z$, so we have three possible cases.
\begin{enumerate}
\item [Case I:] $\nu=0$. In this case, the theorem follows from the corollary to Theorem \ref{rango}.

\item [Case II:] $\mu=0$. In this case, we can suppose that $\nu\neq0$ or $\omega\neq0$.
Otherwise, the theorem is immediate.
For each $\ell\in\Z^p$, we consider the affine transformation $\psi(\ell)$ of $T^3$  given on the covering $\R^3$ by
\begin{equation}
\widetilde{\psi}(\ell)(x,y,z)=(x+\alpha_1(\ell),y+\alpha_2(\ell),\omega(\ell)x+\nu(\ell)y+z+\beta(\ell))
\label{16c}
\end{equation}
Since that we can assume that $\widetilde{\psi}$ is a $\Z^p$-action on $\R^3$, i.e.
\begin{equation}
\widetilde{\psi}(\ell+\ell')=\widetilde{\psi}(\ell)\circ\widetilde{\psi}(\ell')=\widetilde{\psi}(\ell)\circ\widetilde{\psi}(\ell')
\mbox{ for all } \ell,\ell'\in\Z^p
\label{17c}
\end{equation}
and (\ref{17c}) implies that
\begin{enumerate}
\item [i.] $\alpha_1(\ell+\ell')=\alpha_1(\ell)+\alpha_1(\ell')$, for all $\ell,\ell'\in\Z^p$
\item [ii.] $\alpha_2(\ell+\ell')=\alpha_2(\ell)+\alpha_2(\ell')$, for all $\ell,\ell'\in\Z^p$
\item [iii.] $\beta(\ell+\ell')=\omega(\ell)\alpha_1(\ell')+\nu(\ell)\alpha_2(\ell')+\beta(\ell)+\beta(\ell')$, for all $\ell,\ell'\in\Z^p$
\end{enumerate}
Now exchanging $\ell$ and $\ell'$ in item iii above, we have that
 \begin{equation}
 \omega(\ell)\alpha_1(\ell')+\nu(\ell)\alpha_2(\ell')=\omega(\ell')\alpha_1(\ell)+\nu(\ell')\alpha_2(\ell), \mbox{ for all}  \ell,\ell'\in\Z^p.
 \label{18c}
 \end{equation}
 Items i, ii, iii and (\ref{18c}) imply that $\widetilde{\psi}$ in $(\ref{16c})$ is an affine $\Z^p$-action on $\R^3$.
For each $\ell\in\Z^p$ we choose
 \begin{equation}
 \alpha_1(\ell)=a\omega(\ell)\mbox{ and }\alpha_2(\ell)=b\nu(\ell)
 \label{5.1}
 \end{equation}
  and $\beta(\ell)=\frac{1}{2}(\omega(\ell)\alpha_1(\ell)+\nu(\ell)\alpha_2(\ell))+\widehat{\beta}(\ell)$, where $a,b\notin\Q$ and $\widehat{\beta}$ is an homomorphism of $\Z^p$ to $\R$ such that $\widehat{\beta}(\ell)\notin\Z$ for all $\ell\in\Z^p-\{0\}$.
Under these conditions one can see that the affine $\Z^p$-action $\psi$ on $T^3$ induced by $\widetilde{\psi}$ is a free action.

\item [Case III:] $\mu\neq0$ and $\nu\neq0$.
Since ${\bf A}$ is a $\Z^p$-action, $\omega(\ell+\ell')=\mu(\ell)\nu(\ell')+\omega(\ell)+\omega(\ell')$. Now exchanging $\ell$ and $\ell'$, we have that $\mu(\ell)\nu(\ell')=\mu(\ell')\nu(\ell)$.
Hence $\nu=r\mu$ where $r\in\Q$.
For each $\ell\in\Z^p$, we consider the affine transformation $\psi(\ell)$ of $T^3$  given on the covering $\R^3$ by
\begin{equation}
\widetilde{\psi}(\ell)(x,y,z)=(x+\alpha_1(\ell),\mu(\ell)x+y+\alpha_2(\ell),\omega(\ell)x+\nu(\ell)y+z+\beta(\ell))
\label{19c}
\end{equation}
As in the previous case, $\widetilde{\psi}$ satisfies (\ref{17c}), and then we have the following:
\begin{enumerate}
\item [i.] $\alpha_1(\ell+\ell')=\alpha_1(\ell)+\alpha_1(\ell')$, for all $\ell,\ell'\in\Z^p$
\item [ii.] $\alpha_2(\ell+\ell')=\mu(\ell)\alpha_1(\ell')+\alpha_2(\ell)+\alpha_2(\ell')$, for all $\ell,\ell'\in\Z^p$
\item [iii.] $\beta(\ell+\ell')=\omega(\ell)\alpha_1(\ell')+\nu(\ell)\alpha_2(\ell')+\beta(\ell)+\beta(\ell')$, for all $\ell,\ell'\in\Z^p$
\end{enumerate}
Now exchanging $\ell$ and $\ell'$ in ii and iii above, we have that
\begin{enumerate}
\item [a.] $\mu(\ell)\alpha_1(\ell')=\mu(\ell')\alpha_1(\ell)$, for all $\ell,\ell'\in\Z^p$
\item [b.]$\omega(\ell)\alpha_1(\ell')+\nu(\ell)\alpha_2(\ell')=\omega(\ell')\alpha_1(\ell)+\nu(\ell')\alpha_2(\ell)$, for all $\ell,\ell'\in\Z^p$
\end{enumerate}
and now items i, ii, iii, a, and b  imply that $\widetilde{\psi}$ in $(\ref{19c})$ is an affine $\Z^p$-action on $\R^3$.
Let $e_1\in\Z^p$ be such that $\Z^p=\Z e_1\oplus \ker\mu$.
Choose $\alpha_1$ such that $\alpha_1(e_1)\notin \Q $. Then by item a, $\ker\alpha_1=\ker\mu$, and we can define $\alpha_2(\ell)=-\frac{1}{2}\mu(\ell)\alpha_1(\ell)+\frac{\alpha(e_1)}{\nu(e_1)}\omega(\ell)$ for all $\ell\in\Z^p$ and define $\beta(\ell)=-\frac{1}{3}\mu(\ell)\nu(\ell)\alpha_1(\ell)+\frac{1}{2}(\omega(\ell)\alpha_1(\ell)+\nu(\ell)\alpha_2(\ell))+\widehat{\beta}(\ell)$ for all $\ell\in\Z^p$, where $\widehat{\beta}$ is an homomorphism of $\Z^p$ to $\R$ such that $\widehat{\beta}(\ell)\notin\Z$ for all $\ell\in\ker\mu-\{0\}$.
Under these conditions, one can see that the affine $\Z^p$-action $\psi$ on $T^3$ induced by $\widetilde{\psi}$ is a free action.
\end{enumerate}

\subsection{Free minimal affine $\Z^p$-actions on the three-dimensional torus }
We will now analyse which of the actions previously defined can be minimal actions. To this end, we consider the following algebraic characterization of minimal actions, given in ~\cite{L}.
Let $\psi$ be an affine $\Z^p$-action on $T^q$. Denote by $\psi^\ast$ the $\Z^p$-action induce by $\psi$ on the first cohomology group  $H^1(T^q,\Z)=\Z^q$.

Let $\Gamma\subset\Z^p$ be the set of fixed points of $\psi^\ast$.
\begin{equation}
\fix ({\psi^\ast})=\{k\in \Z^{q}|{\bf \psi^\ast}(\ell )k=k \text{ for all }\ell \in \Z^{p}\}.
\end{equation}
We say that $\psi$ satisfies the \textit{irrationality condition on $\Gamma$} if for each $k\in\Gamma-\{0\}$ there is an $\ell\in\Z^p$ such that
\begin{equation}
\langle k,\alpha(\ell)\rangle\notin \Z
\label{5.1i}
\end{equation}
where $\widehat{\psi}(\ell)=\psi_\ast+\alpha(\ell)$ is any lift of $\psi(\ell)$ to the covering $\R^q$.

In Case II, it is easy to see that $\Gamma=\Z^2\times\{0\}$ and that for all $k\in\Gamma$
\begin{equation}
\langle k,\alpha(\ell)\rangle=k_1\alpha_1(\ell)+k_2\alpha_2(\ell)
\label{5.1II}
\end{equation}

If $\{\omega,\nu\}$ is linearly independent over $\Q$ then $\Z^p=\Z e_1\oplus\Z e_2\oplus(\ker\omega\cap\ker\nu)$. Taking the numbers $a,b$ linearly independent over $\Q$ in (\ref{5.1}), we have that for all $k\in\Gamma-\{0\}$ there is $\ell=l_1e_1+l_2e_2$, $l_1,l_2\in\Z-\{0\}$ and the irrationality condition in (\ref{5.1i})  is satisfied.

If $\{\omega,\nu\}$ is linearly dependent over $\Q$ and $(\omega,\nu)\neq(0,0)$, then  $\Z^p=\Z e_1\oplus(\ker\omega\cap\ker\nu)$. Taking the numbers $a,b,c$ and $d$ linearly independent over $\Q$, and  $\alpha_1(\ell)=c\nu(\ell)+a\omega(\ell)$ and $\alpha_2(\ell)=b\nu(\ell)+\omega(\ell)d$, we have that for all $k\in\Gamma-\{0\}$ take $\ell=e_1$, and the irrationality condition in (\ref{5.1i})  is satisfied.

In Case III, it is easy to see that $\Gamma=\Z\times\{0\}\times\{0\}$ and since $\alpha(e_1)\notin\Q$, then for each $k\in\Gamma-\{0\}$ we take $\ell=e_1$, and the irrationality condition in (\ref{5.1i})  is satisfied.

Case I is more interesting. We can see that $\Gamma= \{(k_1,k_2,k_3)\in\Z^3/k_2\mu+k_3\omega=0\}$.
If $\{\mu,\omega\}$ is linearly independent over $\Q$, then $\Gamma=\Z\times\{0\}\times\{0\}$.
On the other side, the affine $\Z^p$-action on $T^3$, $\psi(\ell)={\bf A}(\ell)+{\bf \alpha}(\ell)$,
where ${\bf \alpha}(\ell)=(\alpha_1(\ell),\alpha_2(\ell),\beta(\ell))$ satisfies
$$
\mu(\ell)\alpha_1(\ell')=\mu(\ell')\alpha_1(\ell) \mbox{ and }\omega(\ell)\alpha_1(\ell')=\omega(\ell')\alpha_1(\ell)
$$
then $\ker\mu\subset\ker\alpha_1$ and $\ker\omega\subset\ker\alpha_1$.
Since $\ker\mu +\ker\omega=\Z^p$, it follows that $\alpha_1=0$. Hence the irrationality condition (\ref{5.1i}) is not satisfied.
This is the only case in the torus of dimension 3, where the linear action does not come from minimal affine action.
This shows that, in a certain sense, the theorem \ref{rango} is optimal.

If $\{\mu,\omega\}$ is linearly dependent over $\Q$, then $\Gamma= \{(k_1,k_2,k_3)\in\Z^3/ k_2\mu+k_3\omega=0\}=\Z e_1\oplus \ker L$, where $L$ is the linear map defined by $L(k_2,k_3)=k_2\mu(e_1)+k_3\omega(e_1)$ and $\Z^p=\Z e_1\oplus (\ker\mu\cap \ker\omega)$.
If we take $\alpha_1(e_1)\notin\Q$ and $\{\alpha_2(\ell),\beta(\ell)$ linearly independent over $\Q$ for all $\ell\in\Z^p-\{0\}$, then condition (\ref{5.1i}) is satisfied.
%The following example shows that theorem 5.7.1 can not be extended

\section{Appendix.} In ~\cite{Hi}, the following result is proved.
\begin{thm}\label{hirsch}
Let $G$ be a nilpotent group acting linearly on a finite dimensional vector space; denote the resulting $G$-module by $E$.
If $H^0(G,E)=E^G=0$, then $H^i(G,M)=0$ for all $i\geq0$.
\end{thm}

We give a slightly more general version of this theorem, whose proof is based on theorem \ref{hirsch}.
\begin{thm}\label{extension}
Let $G$ a nilpotent group that acts linearly on the left and right on a finite dimensional vector space $M$.
Suppose that the action on the right is by unipotent transformations of $M$ and that the zero vector is the only vector in $M$ invariant under the action on the left.
Then the cohomology of the resulting $G$-bimodule $M$ is acyclic, i.e. the Hochschild cohomology $HH^*(G,M)=0$
\end{thm}
\noindent{\bf Proof.} The proof proceeds by induction on the dimension of $M$.
If $\dim M=1$, then the right action is trivial, i.e. $m\cdot\ell=m$ for all $\ell\in G$ and for all $m\in M$, hence by [1] $H^*(G,M)=0$.
Let now $\dim M >1$. Define the $G$-module
\begin{equation}
M_1=\{m\in M /m\cdot \ell=m, \mbox{for all } \ell\in G\}.
\end{equation}
Since $G$ is a nilpotent group, it is well known that $M_1\neq 0$, and by the definition of $M_1$, the action on the right is trivial. Moreover, the group of invariants
\begin{equation}
H^0(G,M_1)=M_1^G=0.
\end{equation}
Then, again by [1],
\begin{equation}
H^*(G,M_1)=0.
\label{1a}
\end{equation}
The $G$-bimodule structure on $M$ induces a $G$-bimodule structure on the quotient vector space $M/M_1$, given by

\begin{eqnarray*}
\overline{m}\cdot \ell  &=&\overline{m\cdot \ell },\text{ for }\ell \in G
\text{ and }\overline{m}\in M/M_1 \\
\ell\cdot \overline{m}  &=&\overline{\ell\cdot m },\text{ for }\ell \in G
\text{ and }\overline{m}\in M/M_1
\end{eqnarray*}
it is not hard to show that these operations are well defined.

To apply the induction hypothesis, we will show that the  $G$-bimodule $M/M_1$ satisfies the hypothesis of the theorem. Suppose there is a vector $\overline{m_0}\in M/M_1$ such that
\begin{equation}
\ell\cdot\overline{m_0}=\overline{m_0}, \text{for all } \ell\in G.
\end{equation}
Then $m(\ell)=m_0-\ell\cdot m_0\in M_1$ for all $\ell\in G$.
Hence $m:G\rightarrow M_1$ is a 1-cocycle of the $G$-module $M_1$ and so
it follows from (\ref{1a}) that $m_0\in M_1$. Hence $\overline{m_0}$ is the zero vector.
 On the other hand, let $\ell\in G-\{e\}$ and suppose that $\lambda$ is an eigenvalue of the linear transformation induced by the action on the right $\ell:M/M_1 \rightarrow M/M_!$  $\ell(\overline{m}=\overline{m}\cdot\ell)$, for all $\overline{m}\in M/M_1$. Then there is an $\overline{m_0}\in (M/M_1)-\{\overline{0}\}$ such that $\overline{m_0}\cdot=\lambda\overline{m_0}$
 or equivalently
 \begin{equation}
 m_0\cdot(\ell-I)-m_0(\lambda-1)I\in M_1.
 \label{2a}
  \end{equation}
  Since $\ell$ is unipotent, there is a non-negative integer $d$ such that $m_0\cdot(\ell-I)^d=0$.
If $d=0$, it is trivial that $\lambda=1$.
If $d>0$,
we suppose that $m_0\cdot(\ell-I)^{d-1}\neq 0$. Now multiplying Eq.~(\ref{2a}) by $(\ell-I)^{d-1}$, we obtain
\begin{equation}
-m_0(\lambda-1)(\ell-I)^{d-1}=0
\end{equation}
and then $\lambda=1$, which proves that the linear transformations $\overline{m}\mapsto\overline{m}\cdot\ell$ are unipotent for all $\ell\in G$.
Hence, by induction,
\begin{equation}
HH^*(G,M/M_1)=0
\label{3a}.
\end{equation}
Now, by (\ref{1a}) and (\ref{3a}) and considering the long exact sequence of cohomology groups associated with the exact sequence of $G$-bimodules
\begin{equation}
0\longrightarrow M_1\overset{i}{\longrightarrow } M\overset{\pi}{\longrightarrow} M/M_1\longrightarrow 0
\end{equation}
one can easily show that the cohomology of the $G$-bimodule $M$ is acyclic.
$\hfill\square$
\begin{cor}\label{bimodulo}
Let $ {\bf A}_1$, ${\bf A}_2$ be as in  (\ref{aff1}) and let ${\bf V}$ be as in  (\ref{aff2}).
Then there is a linear map  $W_0: \Q^{q_1}\rightarrow \Q^{q_2}$ such that
${\bf V}(\ell)=W_0{\bf A}_1(\ell)-{\bf A}_2(\ell)W_0$ for all $\ell\in\Z^p$.
\end{cor}
\noindent{\bf Proof.} The $\Z^p$-actions ${\bf A}_1$ and ${\bf A}_2$ induce a $\Z^p$-bimodule structure on the vector space Hom$_\Q(\Q^{q_1},\Q^{q_2})$ given by by composition on the right and by composition on the left respectively. Since that ${\bf A}_1$ is unipotent action and $0$ is the only vector in Hom$_\Q(\Q^{q_1},\Q^{q_2})$ invariant by the action ${\bf A}_2$ then by theorem \ref{extension} $HH^*(\Z^p,\mbox{Hom}_\Q(\Q^{q_1},\Q^{q_2}))=0$ the corollary follows directly from the fact that the first cohomology group $HH^1(\Z^p,\mbox{Hom}_\Q(\Q^{q_1},\Q^{q_2}))=0$.
$\hfill\square$

\section*{}

\bigskip
Richard Urz\'ua Luz
\\ Universidad Cat\'olica del Norte, \\ Casilla 1280, Antofagasta, Chile.
\\ rurzua@ucn.cl
\end{document}